\documentclass{article}
\usepackage[pdftex]{hyperref}
\usepackage{tikz-cd, url}

\hypersetup{
  colorlinks   = true, 
  urlcolor     = gray, 
  linkcolor    = red, 
  citecolor   =  blue
}
\usepackage{amsmath,amsfonts,amsthm,amssymb,graphicx,xcolor}
\usepackage[all,2cell,ps]{xy}

\theoremstyle{plain}
\newtheorem{thm}{Theorem}[section]
\newtheorem{lem}[thm]{Lemma}
\newtheorem{prop}[thm]{Proposition}
\newtheorem{cor}[thm]{Corollary}
\newtheorem{conj}[thm]{Conjecture}

\theoremstyle{definition}

\newtheorem{qtn}[thm]{Question}
\newtheorem{rem}[thm]{Remark}

\newcommand{\R}{\mathbb R}

\newcommand{\C}{\mathbb C}

\newcommand\sH{{\mathcal H}}

\usepackage{mathtools}

\DeclarePairedDelimiter\abs{\lvert}{\rvert}
\DeclarePairedDelimiter\norm{\lVert}{\rVert}

\makeatletter
\let\oldabs\abs
\def\abs{\@ifstar{\oldabs}{\oldabs*}}
\let\oldnorm\norm
\def\norm{\@ifstar{\oldnorm}{\oldnorm*}}
\makeatother

\newenvironment{pf}{\begin{proof}}{\end{proof}}

\title{A rigidity theorem for Einstein $4$-manifolds with sectional curvature of a fixed sign, and its consequences}
 \author{\small{Luca F. Di Cerbo} \\ \scriptsize{University of Florida} \\ \footnotesize{\textsf{ldicerbo@ufl.edu}}}
\date{}

\begin{document}

\maketitle

\begin{abstract}
Any oriented  $4$-dimensional Einstein manifold with semi-definite sectional curvature  (that is, everywhere non-positive or non-negative) satisfies the pointwise inequality
\[
\frac{|s|}{\sqrt{6}}\geq|W^+|+|W^-|,
\]
where $s$, $W^+$ and $W^-$ are respectively the scalar curvature, the self-dual and anti-self-dual Weyl curvatures. We give a complete characterization of closed $4$-dimensional Einstein manifolds with semi-definite sectional curvature realizing the (pointwise) equality case of this inequality. We then present further consequences of this circle of ideas, in particular to the study of the geometry and topology of non-positively curved closed Einstein and K\"ahler-Einstein $4$-manifolds. In the K\"ahler-Einstein case, we obtain a \emph{sharp} Gromov-L\"uck type inequality.
\end{abstract}
\vspace{8cm}
\tableofcontents\quad\\

\vspace{1cm}

\section{Introduction and Main Results}

Given a Riemannian manifold $(M, g)$, we denote by $sec_g$ its sectional curvature, and by $Ric_g$ its Ricci tensor. A Riemannian manifold is said to be Einstein if its Ricci tensor is proportional to the metric
\[
Ric_g=\lambda g,
\] 
where $\lambda\in \R$ is the Einstein constant. The literature concerning the study of Einstein metrics is very vast and continuously growing. We refer the pseudonymous book of Besse \cite{Besse} and the panoramic book of Berger \cite[Chapters 11 \& 12]{Berger} for an introduction to this circle of ideas.

Any oriented $4$-dimensional Einstein metric with semi-definite sectional curvature satisfies the pointwise inequality
\begin{align}\label{semidefinite1}
\frac{|s|}{\sqrt{6}}\geq|W^+|+|W^-|,
\end{align}
where $s$, $W^+$ and $W^-$ are respectively the scalar curvature, the self-dual and anti-self-dual Weyl curvatures. For a proof of this pointwise inequality we refer to Proposition \ref{semidefinite}. This inequality has its roots in the work of LeBrun and Gursky \cite{GLeBrun}, where it is proved for Einstein $4$-manifolds with non-negative sectional curvature and used to show that any Einstein $4$-manifold with non-negative sectional curvature and  and positive intersection form is $\C P^2$ with the Fubini-Study metric.  In this paper, we are mainly concerned about non-positively curved spaces. Section \ref{Proofs} contains a proof of Equation \eqref{semidefinite1} and recalls some of the arguments of Gursky and LeBrun. In Section \ref{Rigidity}, we give a characterization of oriented, closed, Einstein $4$-manifolds realizing the (pointwise) equality case of the inequality in \eqref{semidefinite1}. The second half of the following theorem is the main rigidity theorem contained in this paper. 

\begin{thm}\label{MainT}
Let $(M, g)$ be an oriented Einstein $4$-manifold with Einstein constant $\lambda\in \R$. If the sectional curvature of $g$ is semi-definite, then
\[
\frac{|s_g|}{\sqrt{6}}\geq|W^+|+|W^-|
\]
at every point of $M$. If this inequality is an equality at every point, the universal Riemannian cover $(\widetilde{M}, \tilde{g})$ is completely determined by the Einstein constant:
\begin{equation*}
	(\widetilde{M}, \tilde{g})=\begin{cases} 
	 (S^2\times S^2,  g_{sph}+g_{sph}) &\mbox{if }\quad    \lambda > 0,   \\
	(\R^4, \langle~,~\rangle) & \mbox{if }\quad   \lambda=0,  \\
        (\sH^2_{\R}\times\sH^2_{\R},  g_{hyp}+g_{hyp}) & \mbox{if }\quad \lambda<0.
	\end{cases}
	\end{equation*}
\end{thm}

For a proof of this theorem, we refer to Section \ref{Rigidity}. Finally, we apply this circle of ideas to study the geometry of non-positively curved Einstein $4$-manifolds, see Corollary \ref{finalc}. Remarkably, the result we obtain for non-positively curved K\"ahler-Einstein $4$-manifolds is \emph{sharp}. We refer to Corollary \ref{definitive} and the discussion at the end of Section \ref{geography} for more details. We conclude this study with a question on the geometry of aspherical surfaces of general type that fits naturally in the general theme of Gromov-L\"uck conjecture (see Conjecture \ref{Gromov}).

\noindent\textbf{Acknowledgments}. 
The author thanks Claude LeBrun for a generous correspondence that rekindled this project. He also thanks the two referees for useful bibliographical suggestions and for pertinent comments on the manuscript. He was supported in part by NSF grant DMS-2104662.

\section{Einstein $4$-Manifolds with Semi-Definite Sectional Curvature}\label{Proofs}

In this section, we collect some preliminary results concerning the curvature of Einstein $4$-manifolds with semi-definite sectional curvature. This builds directly on the work of Gursky and LeBrun \cite{GLeBrun}, and we follow closely their arguments and notations. We start with a lemma that is the non-positive sectional curvature analogue of \cite[Lemma 1]{GLeBrun}.

\begin{lem}\label{lemma1}
Let $(M^4, g)$ be an oriented Einstein $4$-manifold. If the sectional curvature of $g$ is non-positive, then
\[
\frac{|s_g|}{\sqrt{6}}\geq |W^+|+|W^-|
\]
at every point of $M$. Moreover, equality holds if and only if
\begin{align}\label{eigenvalues}
\frac{|s_g|}{\sqrt{6}}=\nu_{+}+\nu_{-}, \quad \frac{|W^+|}{\sqrt{6}}=\nu_+=\mu_+, \quad \frac{|W^-|}{\sqrt{6}}=\nu_-=\mu_-.
\end{align}
\end{lem}

\begin{pf}
Let $\mathcal{R}:\Lambda^2\to\Lambda^2$ be the curvature operator. We have that $sec_{g}\leq 0$ if and only if
\begin{equation}\label{LeBrun}
\langle\psi^{+}+\psi^{-}, \mathcal{R}(\psi^{+}+\psi^{-}) \rangle\leq 0
\end{equation}
for all unit-length self-dual $2$-forms $\psi^{+}$ and all all unit-length anti-self-dual $2$-forms $\psi^{-}$. Let $\lambda_{+}\leq \mu_{+}\leq \nu_{+}$ (resp. $\lambda_{-}\leq \mu_{-}\leq \nu_{-}$) be the eigenvalues of $W^{+}$ (resp. $W^{-}$). Recall that $W^{+}$ (resp. $W^{-}$) is a trace-free endomorphsm of $\Lambda^{+}$ (resp. $\Lambda^{-}$). Since $g$ is Einstein, Equation \eqref{LeBrun} implies that
\begin{equation}\label{nonpositive}
\frac{s_g}{6}+\nu^{+}+\nu^{-}\leq 0.
\end{equation}
Next, we observe that
\begin{equation}\label{W+}
|W^{+}|^2=\lambda^2_{+}+\mu^{2}_{+}+\nu^{2}_{+}=2\nu^{2}_{+}-2\lambda_{+}\mu_{+}=2\nu^{2}_{+}+2\nu_{+}\mu_{+}+2\mu^{2}_{+},
\end{equation}
since $\lambda_{+}+\mu_{+}+\nu_{+}=0$. Thus, if $\mu_{+}\geq 0$ in Equation \eqref{W+}, we have $|W^{+}|^2\leq 6\nu^{2}_{+}$ with equality if and only if $\nu_{+}=\mu_{+}$. Similarly, if $\mu_{+}<0$, we have $|W^+|^2=2\nu^2_{+}$. In both cases, we have $|W^+|^2\leq 6\nu^2_{+}$ with equality if and only if $\nu_{+}=\mu_{+}$. Similarly,  we have $|W^-|^2\leq 6\nu^2_{-}$ with equality if and only if $\nu_{-}=\mu_{-}$. If we now use Equation \eqref{nonpositive}, we conclude that
\[
\frac{|s_g|}{6}\geq \nu_{+}+\nu_{-}\geq \frac{1}{\sqrt{6}}\Big(|W^+|+|W^-|\Big).
\]
The proof is complete.
\end{pf}

We can now state the following bound for the Weyl curvatures of Einstein $4$-manifolds with semi-definite sectional curvature.

\begin{prop}\label{semidefinite}
Let $(M^4, g)$ be an oriented Einstein $4$-manifold. If the sectional curvature of $g$ is semi-definite, then
\[
\frac{|s_g|}{\sqrt{6}}\geq |W^+|+|W^-|
\]
at every point of $M$.
\end{prop}
\begin{pf}
If  $sec_{g}\geq 0$, we use Lemma 1 in \cite{GLeBrun}, while if $sec_{g}\leq 0$, we use Lemma \ref{lemma1}.
\end{pf}

\section{A Rigidity Theorem}\label{Rigidity}

We start with a rigidity result for closed, non-flat, non-positively curved, Einstein $4$-manifolds.

\begin{prop}\label{proposition1}
Let $(M^4, g)$ be a closed, oriented, non-flat Einstein $4$-manifold. If the sectional curvature of $g$ is non-positive and
\[
\frac{|s_g|}{\sqrt{6}}=|W^+|+|W^-|
\]
at every point of $M$, then $g$ is locally symmetric and the Riemannian universal cover of $(M, g)$ is the Riemannian product  $(\sH^2_{\R}\times\sH^2_{\R},  g_{hyp}+g_{hyp})$ of two identically-rescaled copies of the hyperbolic plane. In particular, we have $\tau(M)=0$.
\end{prop}

\begin{pf}

In Lemma \ref{lemma1}, we showed that if the equality
\[
\frac{|s_g|}{\sqrt{6}}=|W^+|+|W^-|
\]
holds then
\begin{align}\label{eigenvalues}
\frac{|W^+|}{\sqrt{6}}=\nu_+=\mu_+, \quad \frac{|W^-|}{\sqrt{6}}=\nu_-=\mu_-.
\end{align}
Next, we prove that the metric $g$ cannot be self-dual or anti-self-dual. By contradiction assume $|W^-|=0$. Since $g$ is assumed to be non-flat, we must have that $W^+\neq 0$ at each point of $M$.  Moreover, because of Equation \eqref{eigenvalues}, $W^+$ has less than three distinct eigenvalues at each point of $M$. Thus, by \cite[Proposition 5]{Der83} we have that either $(M, g)$ or a double Riemannian cover of it is a globally conformally K\"ahler $4$-manifold and $W^+\neq 0$ at every point of $M$. Since any Riemannian cover of $(M, g)$ has the same universal Riemannian cover as $(M, g)$, we can, without loss of generality, assume that $(M, g)$ is globally conformally K\"ahler.  Hadamard's theorem (\emph{e.g.,} \cite[Theorem 6.2.2]{Petersen}) now tells us that the topological universal cover of $M$ is diffeomorphic to $\R^4$. In other words, $M$ is an aspherical complex surface. Thus, by \cite[Lemma 2]{ADC23} $(M, J)$ is a minimal complex surface. Moreover, Chern-Gauss-Bonnet and Thom-Hirzebruch integral formulas (\emph{e.g.,} \cite[Section 2]{LeB99})  imply that
\[
\chi(M)>0, \quad \tau(M)>0 \quad \Rightarrow \quad c^2_{1}(M, J)=2\chi(M)+3\tau(M)>0.
\]
By \cite[Lemma 3]{ADC23}, we have that $(M, J)$ is a minimal surface of general type with ample canonical line bundle. In particular, $(M, J)$ satisfies the celebrated Bogomolov-Miyaoka-Yau inequality
\begin{align}\label{BMY}
\chi(M)\geq 3\tau(M)(>0),
\end{align}
for which we refer to classical book \cite[Chapter VII, Section 4]{BHPV04}. On the other hand, if we combine the Chern-Gauss-Bonnet and Thom-Hirzebruch integral formulas, with the fact that $g$ is Einstein and it satisfies
\[
\frac{|s_g|}{\sqrt{6}}=|W^+|>0, \quad |W^-|=0,
\]
we have
\begin{align}\label{CGBHT}
\chi(M)=\frac{1}{8\pi^2}\int_{M}|W^+|^2+\frac{s^2_g}{24}d\mu_g=\frac{15}{8}\frac{1}{12\pi^2}\int_{M}|W^+|^2d\mu_g=\frac{15}{8}\tau(M).
\end{align}
Equation \eqref{CGBHT} contradicts the bound given in Equation \eqref{BMY}, and the argument is complete. Similarly, if we assume $|W^+|=0$, we can change orientation and argue as before. Concluding, we must have $W^+\neq 0$ and $W^-\neq 0$ at every point of $M$ (\cite[Proposition 5, part iv]{Der83}). 

Next, with respect to the first orientation, $(M^4, g)$ is globally conformally K\"ahler. Let us denote by $(M, J, \bar{g})$ such K\"ahler surface. Since every $4$-dimensional conformally Einstein metric is Bach-flat (see for more details \cite[Section 2]{LeB20}), we have that $\bar{g}$ must have vanishing Bach tensor. Thus, $\bar{g}$ is a critical point of the Weyl functional
\[
\mathcal{W}(g):=\int_{M}|W|^2_g d\mu_g=-12\pi^2\tau(M)+2\int_{M}|W^+|^2_g d\mu_{g}.
\]
Since for a K\"ahler metric we have
\[
|W^+|^2_{\bar{g}}=\frac{s^2_{\bar{g}}}{24},
\]
we deduce that $\bar{g}$ is a critical point of the Calabi functional
\[
\mathcal{C}(\omega):=\int_{M}s^2_{\omega} d\mu_{\omega}
\]
for $\omega$ in the K\"ahler class $\Omega:=[\omega_{\bar{g}}]$ (but more generally also on the whole space of K\"ahler metrics). In other words, $\omega_{\bar{g}}$ is an extremal K\"ahler metric in the sense of Calabi. This implies that $X:=J\nabla s_{\bar{g}}$ is a Killing field of $\bar{g}$ (this also follows directly from the local argument in \cite[Proposition 4]{Der83}), and also that 
\begin{align}\label{uno}
X(s_{\bar{g}})=0.
\end{align}
Importantly, we have that $X$ is a Killing field for $g$ as well. To see this we argue as follows. First recall that following \cite[Proposition 5]{Der83}, we have that 
\[
g=\bar{g}s^{-2}_{\bar{g}}.
\]
Since $X$ is Killing for $\bar{g}$, we have
\begin{align}\label{due}
X_{j, k}+X_{k, j}=\bar{g}_{js}\partial_{k}X^s+\bar{g}_{ks}\partial_{j}X^s+X^s\partial_s\bar{g}_{jk}=0.
\end{align}
This implies that 
\begin{align}\notag
g_{js}\partial_{k}X^s+g_{ks}\partial_{j}X^s+X^s\partial_s g_{jk}=&s^{-2}_{\bar{g}}\big(\bar{g}_{js}\partial_{k}X^s+\bar{g}_{ks}\partial_{j}X^s+X^s\partial_s\bar{g}_{jk}\big)\\ \notag
&-2s^{-3}_{\bar{g}}X(s_{\bar{g}})\bar{g}_{jk}=0
\end{align}
because of Equations \eqref{uno} and \eqref{due}. Thus, we have $X_{j, k}+X_{k, j}=0$ with respect to the metric $g$ as well. This shows that $X$ is a Killing for $(M, g)$.
Next, since $(M^4, g)$ has negative Ricci curvature and it is compact with no boundary, Bochner's formula and an integration by part give that we cannot have any non-trivial Killing fields for $g$. Thus, $s_{\bar{g}}$ must be constant and then $g$ coincides with $\bar{g}$ up to scaling. Concluding $(M, g)$ is K\"ahler-Einstein and the self-dual part of the Weyl tensor is parallel. Now, the same argument holds with the reverse orientation, so that the Weyl tensor is parallel and $g$ must then be locally symmetric. Because of the sign of the Einstein constant, we have that the Riemannian universal cover of $(M, g)$ is $(\sH^2_{\R}\times\sH^2_{\R},  g_{hyp}+g_{hyp})$. Since $|W^+|=|W^-|$, we deduce that $\tau(M)=0$.
\end{pf}

We also have the non-negative version of this rigidity result. The argument is somewhat different, but it still involves a mixture of local and global geometry.

\begin{prop}\label{proposition2}
Let $(M^4, g)$ be a closed, oriented, non-flat Einstein $4$-manifold. If the sectional curvature of $g$ is non-negative and
\[
\frac{s_g}{\sqrt{6}}=|W^+|+|W^-|
\]
at every point of $M$, then $g$ is locally symmetric and the Riemannian universal cover of $(M, g)$ is the Riemannian product $(S^2\times S^2,  g_{sph}+g_{sph})$ of two identically-rescaled copies of the standard $2$-sphere. In particular, we have $\tau(M)=0$.
\end{prop}

\begin{pf}
By \cite[Lemma 1]{GLeBrun}, we have that
\[
\frac{s_g}{\sqrt{6}}=|W^+|+|W^-|
\]
implies 
\[
\frac{|W^+|}{\sqrt{6}}=|\lambda_+|=|\mu_+|, \quad \frac{|W^-|}{\sqrt{6}}=|\lambda_-|=|\mu_-|,
\]
at every point of $M$. We start by selecting the orientation on $M$ so that $\tau(M)\geq 0$, and we rule out the case $W^-=0$. First, by Equation \eqref{CGBHT} we observe that
\[
\chi(M)=\frac{15}{8}\tau(M).
\]
Bourguignon formula \cite{Bou81} for harmonic $2$-forms with on self-dual $4$-manifolds with $s>0$ implies that $b^-_{2}=0$. Now because of the Myers' theorem (\emph{e.g.,} \cite[Chapter 9]{Petersen}), we have that $b_1=0$. Thus, $\chi(M)=2+\tau(M)$ which implies $\tau(M)=\frac{16}{7}$. This is a contradiction as the signature is an integer.

By reiterating the same argument with the opposite orientation, we conclude once again because of  \cite[Proposition 5]{Der83} that $W^+\neq 0$ and $W^-\neq 0$ at every point of $M$. Moreover, we have that both $W^+$ and $W^-$ have less than three distinct eigenvalues at every point of $M$. Without loss of generality, we now assume $M$ to be simply connected. By \cite[Proposition 5]{Der83}, if $g$ is not K\"ahler-Einstein then $g$ is globally conformally K\"ahler. By \cite[Theorem A]{LeB12},  $(M, g)$ admits an orientation compatible complex structure $J$ that makes $(M, J)$ a \emph{del Pezzo} surface with $g$ conformally K\"ahler. Since we can repeat the same argument with the opposite orientation, if $(M, g)$ is K\"ahler-Einstein then it has to be locally symmetric and then $(M, g)$ is isometric to  $(S^2\times S^2,  g_{sph}+g_{sph})$.

If $g$ is not K\"ahler-Einstein, then again by \cite[Theorem A]{LeB12}, it has to be either the Page metric \cite{Page} on $\C P^2\#\overline{\C P^2}$  or the Chen-LeBrun-Weber metric \cite{CLW08} on  $\C P^2\#2\overline{\C P^2}$. Now, the Page metric has some negative curvature as it can be shown by a \emph{Maple} computation, see \cite[page 31]{KK18}. Finally, since $c^2_{1}$ is \emph{not} even for $\C P^2\#2\overline{\C P^2}$, we can rule out the case of the Chen-LeBrun-Weber metric as the underlying smooth manifold cannot support two distinct complex structures compatible with the distinct orientations. This observation is due to Beauville \cite{Bea85}, and we recall its proof here for completeness sake. Let $(M, J)$ be a complex surface and assume the underlying smooth $4$-manifolds with the opposite orientation admits compatible complex structure, say $(\overline{M}, \bar{J})$. By the Todd genus formula, we have
\[
c^2_{1}(\overline{M}, \bar{J})\equiv c^2_{1}(M, J)\mod 12.
\]
Since we also that
\[
c^2_{1}(\overline{M}, \bar{J})-c^2_{1}(M, J)=-6\tau(M),
\]
we conclude that $\tau(M)$ is even and then 
\[
c^2_{1}(M, J)=2\chi(M)+3\tau(M)\equiv 0 \mod 2.
\]
\end{pf}

We can now prove the main rigidity result stated in the introduction

\begin{pf}[Proof of Theorem~\ref{MainT}]
The statement is now a direct consequence of Propositions \ref{semidefinite}, \ref{proposition1} and \ref{proposition2}.
\end{pf}

\section{Geometry of Non-Positively Curved Einstein $4$-Manifolds}\label{geography}

As discussed by Gromov in \cite[Section 8]{Gro93} and L\"uck \cite[Theorem 5.1]{Luck}, the famous conjecture of I. Singer on the vanishing of $L^2$-Betti numbers of aspherical manifolds implies a constraint on the geometry and topology of such manifolds in dimension four. In more details, the Singer conjecture implies a \emph{refined} version of the celebrated problem of Hopf on the sign of the Euler characteristic of aspherical, closed, $4$-manifolds. 

\begin{conj}[Gromov-L\"uck Inequality]\label{Gromov}
	If $X$ is a closed, oriented, aspherical $4$-manifold, then:
	\begin{equation*}
	\chi(X)\geq |\tau(X)|,
	\end{equation*}
where $\tau(X)$ is the signature of $X$.
\end{conj}

\begin{rem}
Surprisingly, the conjecture of Singer is currently open for complex surfaces, and there is a large variety of aspherical complex surfaces, \textit{e.g.}, most smooth toroidal compactifications of ball quotient surfaces \cite[Theorem A]{DiC12}. That said, it follows from \cite{ADL23} and \cite{DL24} that the conjecture of Singer holds true for aspherical complex surfaces with residually finite fundamental group.
\end{rem}

The inequality in Conjecture \ref{Gromov} is very coarse, and indeed there is no conjectural statement regarding the equality case. Interestingly, it seems to be currently unknown whether aspherical $4$-manifolds with $\chi=|\tau|$, $\tau\neq 0$ exist. Also, frustratingly enough, we do not know any example of aspherical $4$-manifolds with $3|\tau|>\chi>0$. For more on this circle of ideas we refer to \cite[Section 1]{DG23}, \cite[Section 1]{ADL23}, and \cite[Section 3]{DP23}.

In this section, we observe that our study (\emph{cf}. Sections \ref{Proofs} and \ref{Rigidity}) can shed some light on Conjecture \ref{Gromov} for non-positively curved Einstein $4$-manifolds. Moreover, we provide an \emph{improved} and \emph{sharp} version of Gromov-L\"uck inequality for non-positively curved K\"ahler-Einstein surfaces.  

We begin by a result that was originally hinted in the beautiful survey paper \cite[Remark 10.1]{LeB99}.

\begin{cor}\label{finalc}
Let $(M^4, g)$ be a closed, oriented, non-flat Einstein $4$-manifold. If the sectional curvature of $g$ is non-positive, then
\[
\chi(M)>\frac{15}{8}|\tau(M)|.
\]
\end{cor}

\begin{pf}
Take the orientation such that $\tau(M)\geq 0$. If $\tau(M)=0$, then there is nothing to prove as the Euler characteristic of a non-flat Einstein metric is strictly positive. By the same curvature computations as in Proposition \ref{proposition1}, we immediately get that
\begin{align}\label{preliminary}
\chi(M)\geq \frac{15}{8}\tau(M),
\end{align}
with equality if and only if $|W^+|=\frac{|s_g|}{\sqrt{6}}$ and $|W^-|=0$. Up to possibly a double cover (see once again  \cite[Proposition 5]{Der83}), $M$ is now also globally conformally K\"ahler. Since $c^2_1>0$, $(M, J)$ is a minimal surface of general type (see \cite[Lemma 2]{ADC23}). This implies that the Bogomolov-Miyaoka-Yau inequality $\chi(M)\geq 3\tau(M)>0$ must hold. Concluding, the inequality in Equation \eqref{preliminary} must be strict. 
\end{pf}

We conclude with a sharp inequality for non-positively curved K\"ahler-Einstein surfaces.

\begin{cor}\label{definitive}
Let $(M, J, g)$ be a K\"ahler-Einstein surface of non-positive sectional curvature. We then have $\tau(M)\geq 0$ with equality if and only if the universal Riemannian cover $(\widetilde{M}, \tilde{g})$ is either $(\C^2, |dz_1|^2+|dz_2|^2)$ or $(\sH^2_{\R}\times\sH^2_{\R}, g_{hyp}+g_{hyp})$. Moreover, if $(M, g)$ is non-flat then $(M, J)$ is a minimal surface of general type satisfying 
\[
\chi(M)>0, \quad \textrm{ } \quad \chi(M)\geq 3\tau(M),
\]
with equality $\chi(M)=3\tau(M)$ if and only of $M$ is a complex-hyperbolic surface (i.e., its Riemannian universal cover is the unit ball in $\C^2$ equipped with the symmetric Bergman metric say $(\sH^2_{\C}, g_{Ber})$). 
\end{cor}
\begin{pf}
For a K\"ahler surface, we always have
\[
|W^+|^2=\frac{s^2}{24}\quad \Rightarrow \quad |W^+|=\frac{|s|}{2\sqrt{6}}.
\]
By Lemma \ref{lemma1}, we have the pointwise inequality
\[
|W^-|\leq\frac{|s|}{2\sqrt{6}}\quad \Rightarrow \quad |W^-|^2\leq\frac{s^2}{24}.
\]
We conclude that the integrand in the Hirzebruch-Thom formula for the signature is non-negative, so that we obtain
\[
\tau(M)=\frac{1}{12\pi^2}\int_{M}|W^+|^2-|W^-|^2d\mu_g\geq 0.
\]
Concerning the equality case $\tau(M)=0$, we must have
\[
|W^-|=\frac{|s|}{2\sqrt{6}}=|W^+| \quad \Rightarrow \quad \nabla\mathcal{R}=0,
\]
so that $g$ is locally symmetric (\emph{cf.} \cite[Proposition 9]{Der83}). This implies that the Riemannian universal cover of $(M, g)$ is either $(\C^2, |dz_1|^2+|dz_2|^2)$ or $(\sH^2_{\R}\times\sH^2_{\R}, g_{hyp}+g_{hyp})$.

Next, we assume $(M, g)$ to be non-flat. By Chern-Gauss-Bonnet, we know that $\chi(M)>0$ and by \cite[Lemma 2]{ADC23} we know that $(M, J)$ is a minimal surface of general type. For such surfaces, the Bogomolov-Miyaoka-Yau inequality gives that
\[
\chi(M) \geq 3\tau(M),
\]
with equality of and only if $M$ is a ball quotient.
\end{pf}

\begin{rem}
Many examples of non-locally symmetric, negatively curved, K\"ahler-Einstein manifolds appeared in the recent preprint of Guenancia and H\"amenstadt \cite{GH25}. It may be the case that such objects exist in great numbers, maybe even on any of the previously constructed examples of negatively curved K\"ahler surfaces. We refer to \cite{GH25} for more background on this problem and for countable many examples in \emph{any} complex dimension greater than or equal to two.
\end{rem}

As discussed in the recent paper \cite[Section 3]{DP23}, we currently do not know examples of aspherical complex surfaces of general type with $\tau<0$. Corollary  \ref{definitive} provides some support for the non-existence of such creatures in the non-positively curved setting. That said, the result crucially relies on the K\"ahler-Einstein assumption! Can this assumption be removed?  At the time of writing, it seems we do not know any counterexample to the following question.

\begin{qtn}
Let $(M, J)$ be an aspherical surface of general type. Is it true
\[
\chi(M)\geq 3\tau(M)\geq 0
\]
with $\tau(M)=0$ if and only if $M$ is a quotient $\sH^2_{\R}\times \sH^2_{\R}$?
\end{qtn}

\section{Data Avalability}

There is no data associated to this work.

\section{Conflict of Interest}

The author states that there is no conflict of interest.

%
%
%
%
%
%
%

\end{document}